\documentclass[12pt]{article}
\usepackage{amssymb,amsmath,graphicx,ucs,hyperref}
\usepackage[utf8x]{inputenc}

\newtheorem{definition}{Definition}
\newtheorem{remark}{Remark}

\newtheorem{conjecture}{Conjecture}
\title{How not to scramble a Rubik's cube}
\author{Thomas Fernique}

\begin{document}

\maketitle

\begin{abstract}
We model the scrambling of a Rubik's cube by a Markov chain and introduce a stopping time $T$ which is a quite natural candidate to be a strong uniform time.
This may pave the way for estimating the number of moves required to scramble a cube.
Unfortunately, we show that $T$ is not this is not strongly uniform.
\end{abstract}

Playing with a Rubik's cube involves starting with a ``random'' configuration and, in the shortest possible time, achieving the ``sorted'' configuration by turning the faces.
It is known that 20 moves (counting quarter turns and half turns) are always sufficient and sometimes necessary to solve a cube \cite{RKDD14} (although, in practice, players make more moves because they prefer to perform short sequences of moves whose effect on the distance to the arranged configuration is easier to understand).

But how do we obtain the ``random'' starting configuration that the player will have to solve?
The natural idea is to make random moves ``long enough.''
But can we quantify exactly how long?
We suspect that at least 20 are needed, because otherwise certain initial configurations will never be reached.
Actually, it is known that at least 26 are needed \cite{QRY25}.
To the best of our knowledge, this is the only rigorous result on this subject.
In particular, no upper bound seems to be known.


\begin{remark}
In competition, the initial configuration is not obtained by randomly mixing a cube.
Instead, all the possible configurations are numbered, a program draws a random number and calculates a sequence of moves to reach the corresponding configuration, which a preparer will perform on the cube.
\end{remark}

Formally\footnote{The reader is referred to \cite{LPW06} for classical detailed definitions and results.}, the scrambling is modeled by a {\bf Markov chain} $X_t$, and ``how much'' a cube is scrambled starting from a configuration $x_0$ is measured by the {\bf distance to uniform distribution}:
where $\Omega$ is the set of all the Rubik's cube configurations, $|\Omega$ its size and $\mathbb{P}(X_t=x|X_0=x_0)$ is the probability that the configuration $x_0$ has been transformed into configuration $x$ after $t$ step of the Markov chain.

\begin{remark}
The choice of $x_0$ does not play any role here because of the symmetry of $\Omega$.
The configuration where the small cubes on each face are all the same color is indeed special for us only by its aesthetics.
\end{remark}

For suitable Markov chains (namely ergodic and with a uniform stationnary distribution), $d(t)$ tends towards $0$.
The speed of convergence is then usually measured by the {\bf mixing time}:
$$
\tau:=\min\{t~|~d(t)\leq 1/4\}.
$$
The threshold $1/4$ is arbitrary but classic.
This can be viewed as the half-life: making $\tau$ further moves will halve the distance to the uniform distribution.

In \cite{QRY25}, the authors propose the Markov chain which, at each step, chooses one of the six faces of the cube uniformly at random and rotates it one, two, or three quarters of a turn uniformly at random.
Then, they prove that the mixing time is at least $26$.

Of course, the simplest way to calculate the mixing time is to explicitly calculate $d(t)$ for $t=1,2,3\ldots$ until it falls below $1/4$.
To do this, it is sufficient to maintain, for each configuration $x\in \Omega$, the probability that $X_t=x$.
Unfortunately, the size of $\Omega$ is far too large for this to be feasible.

\begin{remark}
In \cite{QRY25}, the authors introduce a projection from $\Omega$ to a smaller set and compute the mixing time in this smaller set, which is a lower bound on the mixing time over $\Omega$.
\end{remark}

\begin{remark}
The $2\times 2\times 2$ variant of the Rubik's cube (Pocket cube) has much less configurations.
An array of size $|\Omega|$ can easily be maintained on a computer and yields $d(t)$ for every $t$.
In particular, this yields $\tau=19$ for the analogue of the Markov chain defined in \cite{QRY25} (to be compared with the $11$ moves which are always sufficient and sometimes necessary to solve such a cube).
This also yields a lower bound for the mixing time of the usual Rubik's cube, because scrambling it properly requires at least scrambling its corners properly, which is equivalent to scrambling the Pocket cube properly.
\end{remark}

Now, consider the same Markov chain as in \cite{QRY25}.
W.l.o.g, the center of a Rubik's cube as well as the central face of each face can be considered to be fixed.
In other words, scrambling a Rubik's cube amounts to move around $20$ of its small cubes.
We also set an $(X, Y, Z)$ coordinate system whose axes are parallel to the edges of the cube.
Here is the main definition of this paper:

\begin{definition}
A pair of small cubes is said to be {\bf $X$-unlinked} at step $t$ if, at some previous step of the Markov chain, a face perpendicular to the $X$-axis which contains exactly one of the two small cubes has been rotated.
\end{definition}

This definition extends to the axes $Y$ and $Z$, and the pair is said to be {\bf unlinked} if it is $X$-, $Y$- and $Z$-unlinked.
The underlying idea is that small cubes have “forgotten” their initial positions when they are unlinked, and therefore any configuration can be obtained with the same probability when all pairs are unlinked.
Formally, let $T$ be the random variable
$$
T=\min\{t~|~\textrm{in $X_t$, every pair of small cubes is unlinked}\}.
$$
It is a {\bf stopping time} (it depends only on the past of the Markov chain), and one might wonder if this is not a {\bf strong uniform time}, that is, at time $t\geq T$ every configuration has the same probability to appear:
$$
\mathbb{P}(X_t=x|t\geq T)=\frac{1}{|\Omega|}.
$$

Strong uniform times have been introduced in \cite{AD86} in order to bound the mixing time of a deck of cards via the {\bf riffle shuffle}, where the deck is several times split into two halves and interleaved.
Namely, the authors proved the time at which each pair of cards have been separated at least once by the split of the deck is a strong uniform stationary time.
Indeed, the show that any two cards are equally likely to be placed one below the other or vice versa, i.e., any configuration (permutation) of the deck of cards has the same probability of having been obtained.
The connection with the mixing time is then made (still in \cite{AD86}) via the following inequality, valid for any strong uniform time $T$:
$$
d(t)\leq\mathbb{P}(t<T).
$$

The strong uniform stationary time here proposed for scrambling the Rubik's cube is obviously inspired by the above one for shuffling cards.
Since rotating a face moves corners onto corners, we focus on corners.
We would like to prove:

\begin{conjecture}
Fix an arbitrary order on the positions of the corners of a Rubik's cube.
Then, when every pair of corners of $X_t$ are unlinked, any two corners are equally likely to be placed one before the other or vice versa.
\end{conjecture}

Unfortunately, the above conjecture does not hold.
We prove it as follows.
Fix an arbitrary order on the positions of the corners and consider a pair $(a_0,b_0)$ of corners.
We will compute the probability that, at step $t=T$ (our conjectured strong uniform time), theses corners are still in the same order or have been exchanged.
For that:
\begin{enumerate}
\item We build the directed graph $G$ whose vertices are the pairs $(a,b)$ of corner positions, with a directed edge from $(a,b)$ to $(a',b')$ if the Markov chain can bring in one step the corners in position $a$ and $b$ to, respectively, $a'$ and $b'$.
\item Starting from $(a_0,b_0)$, we recursively follows the directed edges and add to every vertex $(a,b)$ a triple $(x,y,z)$ of boolean which specify whether the pair $(a,b)$ has been $X$-, $Y$- or $Z$-unlinked, depending on the face move which bring us to this pair.
\item We start the above recursion by adding $x=y=z=False$ to $(a_0,b_0)$ and stop it as soon as we reach a vertex with $x=y=z=True$.
In this later case we add a loop on the vertex to make it a sink.
\item We put probability $1$ on the vertex $(a_0,b_0,x=y=z=False)$ and we iteratively spread the mass along edges.
The mass is equally distributed among the edges originating from the same vertex.
The mass gradually concentrates in the sinks, and we stop the iteration when the mass outside the sinks is below a certain small threshold $\varepsilon$.
\end{enumerate}
After that, the sum of the mass over the vertices $(a,b,x=y=z=True)$ such that $a<b$ gives the probability within $\varepsilon$ that $a_0$ is before $b_0$ at time $t=T$.
The numerical values obtained are close to $1/2$ but not equal to it: the difference from 0.5 exceeds 0.026 for certain pairs $(a_0,b_0)$.

\begin{remark}
This does not completely rule out the possibility that $T$ is a strong uniform time.
One could indeed imagine that when the very last pair of corners becomes unlinked, the other pairs (which have continued to be moved while waiting for the last one) are fairly ordered.
That seems however very unlikely.
Moreover, we shall also consider the permutation on the remaining small cubes (edges) and the orientations of all the small cubes.
\end{remark}

\begin{remark}
The fact that the probability of a corner being placed before another rather than after is numerically close to $1/2$ is not very surprising.
Simulations show that the expected time to unlink all pairs of corners is approximately $27$.
This is true for both the Rubik's cube and the Pocket cube, since only the corners are considered.
However, the mixing time for the Pocket Cube, explicitly calculated, is $19$.
So at $t=T$ we are already very close to the uniform distribution\ldots

Let us also mention that simulations show that the expected time to unlink all pairs of cubes (not just the corners) in the Rubik's Cube is approximately $41$.
Furthermore, these same simulations combined with the inequality $d(t)\leq\mathbb{P}(t<T)$ would suggest a mixing time of at most $31$ for the Pocket Cube (compared to the actual $19$) and at most $46$ for the Rubik's Cube (compared to the lower bound of $26$ in \cite {QRY25}).
\end{remark}

\begin{remark}
As for shuffling cards, a clever method (based on ascending sequences) was later developed in \cite{BD92} to calculate $d(t)$ exactly, even though the number of possible configurations of a 52-card deck is even greater than that of a Rubik's cube.
Perhaps we should not lose all hope of calculating $d(t)$ exactly for the Rubik's cube\ldots
\end{remark}

\appendix
\section{Code}

\begin{verbatim}
# Faces are denoted by B,D,F,L,R,U
# for (Back, Down, Front, Left, Right, Up)
# Each corner is denoted by the three faces it belongs to.
# They are arbitrarily numbered as follows:
(DFR,FRU,BRU,BDR,FLU,DFL,BLU,BDL)=(0,1,2,3,4,5,6,7)

# each line of M give the corner each corner is mapped to
# by a face quarter turn. For example, M[1][BDR]=BDL
# because a quarter turn of D maps BDR onto BDL
M=[
# DFR, FRU, BRU, BDR, FLU, DFL, BLU, BDL
[ DFR, FRU, BLU, BRU, FLU, DFL, BDL, BDR], #B
[ BDR, FRU, BRU, BDL, FLU, DFR, BLU, DFL], #D
[ DFL, DFR, BRU, BDR, FRU, FLU, BLU, BDL], #F
[ DFR, FRU, BRU, BDR, DFL, BDL, FLU, BLU], #L
[ FRU, BRU, BDR, DFR, FLU, DFL, BLU, BDL], #R
[ DFR, FLU, FRU, BDR, BLU, DFL, BRU, BDL], #U
]

# given two small cubes a, b, boolean x, y, z, face f
# update x, y, z depending on whether rotating f
# has unlink a, b along x, y, z
def update(a,b,x,y,z,f):
    BF=[FLU,FRU,DFR,DFL] # the corners moved by F
    if f in [0,2] and ((a in BF)^^(b in BF)): x=True
    LR=[FLU,BLU,BDL,DFL] # the corners moved by L
    if f in [3,4] and ((a in LR)^^(b in LR)): y=True
    DU=[FLU,BLU,BRU,FRU] # the corners moved by U
    if f in [1,5] and ((a in DU)^^(b in DU)): z=True
    return (a,b,x,y,z)

# build the graph G defined in the text.
def build(a, b, x=False, y=False, z=False):
    # already in G -> cycle -> stops
    if (a,b,x,y,z) in G.keys():
        return
    # otherwise add the vertex to G
    G.update({(a,b,x,y,z):set()})
    # completly unlinked -> make a self-loop and stops
    if x and y and z:
        G[(a,b,x,y,z)].add((a,b,x,y,z))
        return
    # add the vertices than can be reached
    # then recursively explore them
    for face in range(6):
        (a2,b2,x2,y2,z2)=update(a,b,x,y,z,face)
        for q in range(4):
            (a2,b2)=(M[face][a2],M[face][b2])
            G[(a,b,x,y,z)].add((a2,b2,x2,y2,z2))
            build(a2,b2,x2,y2,z2)

# Probability that a0 is before b0
# when the conjectured strong stationnary time is reached
def explore(a0, b0, nsteps=100):
    global G
    G={}
    build(a0,b0)
    W={(a0,b0,False,False,False):1} # initial distribution
    for _ in range(nsteps):
        W2={v:0 for v in W.keys()} # W2 will be the new W
        for u in W.keys(): # u gives weight to neighbors
            for v in G[u]:
                if v in W2.keys():
                    W2[v]+=RR(W[u]/len(G[u]))
                else:
                    W2.update({v:RR(W[u]/len(G[u]))})
        W=copy(W2)
    # probability z that a0 and b0 are totally unlinked
    z=sum([W[i] for i in W if i[2] and i[3] and i[4]])
    print("%.10f of the distribution in the sinks"%(100*z))
    # probability q that a0 is before b0 (should be 0.5)
    q=sum([W[i] for i in W if i[2] and i[3] and i[4]
                            and i[0]<i[1]])
    # bound q depending on z
    return q+(1-z) if q<0.5 else (1-q)-(1-z)
    if q<0.5:
        print("P[%d<%d]<%f"%(a0,b0,q+(1-z)))
    else:
        print("P[%d<%d]<%f"%(b0,a0,(1-q)-(1-z)))
\end{verbatim}

\bibliographystyle{alpha}
\bibliography{rubik}

\end{document}